\documentclass[12pt, reqno]{amsart}
\allowdisplaybreaks
\begin{document}
\setcounter{page}{1}

\title[\hfilneg \hfil Ostrowski type inequalities time scale ]
{Some Ostrowski type inequalities for double integrals  On Time Scales}

\author[Deepak B. Pachpatte\hfil \hfilneg]
{Deepak B. Pachpatte}

\address{Deepak B. Pachpatte \newline
 Department of Mathematics,
 Dr. Babasaheb Ambedkar Marathwada University, Aurangabad,
 Maharashtra 431004, India}
\email{pachpatte@gmail.com}


\begin{abstract}
  The main objective of this paper is to study some Ostrowski and Trapezoid type inequalities for double integrals  on Time Scales.
	Some other interesting inequalities are also given.
 
\end{abstract}
\maketitle

\section{Introduction}
     In  year 1988 the German mathematicain in his Ph.D dissertation has initiated the study of time scales calculus  which unifies the theory of both differential and difference calculus \cite{HIG}. Dynamical equations and inequality's can be used studying various properties  and model many phenomena in economics \cite{Ati}, biological systems \cite{Th} and various systems  in  neural network \cite{Li}.
   
   In \cite{Boh3} Bohner and Matthews have given the Ostrowski inequality and Montgomery identity on time scales. Some results on Ostrowski and Gruss inequality were obtained by N. Ahmad, W. Liu and others \cite{Aha1, Liu1, Tun}. Recently in \cite{Ana, Jia, Liu2, Liu4, Xu} authors have obtained some new  Ostrowski type inequalities. Weighted Ostrowski and Trapezoid inequalities on time scales are  obtained by W. Liu and others  in \cite{Liu5, Liu6, Liu7}. In \cite{Sar} M. Sarikya have studied some weighted Ostrowski and Chebsev type inequalities on time scales. Motivated by the results in the above paper we obtain some Ostrowski  and Trapezoid type inequalities for double integrals on time scales. 
	
	   In what follows the time scale  $\mathbb{T}$ is a nonempty closed subset of $\mathbb{R}$. Let $t \in \mathbb{T}$ the mapping $\sigma ,\rho :\mathbb{T} \to \mathbb{T}$ are defined as $\sigma \left( t \right) = \inf \left\{ {s \in \mathbb{T}:s > t} \right\}$ and $\rho \left( t \right) = \sup \left\{ {s \in \mathbb{T}:s < t} \right\}$ are called the forward and backward jump operators respectively.
		
		We say that $f:\mathbb{T} \to \mathbb{R}$ is rd-continuous provided $f$ is continuous at each right-dense point of $\mathbb{T}$ and has a finite left sided limit at each left dense point of $\mathbb{T}$. $C_{rd}$ denotes the set of rd-continuous function defined on $\mathbb{T}$. Let $\mathbb{T}_1$ and $\mathbb{T}_2$ be two time scales with at least two points and consider the time scales intervals $\overline {\mathbb{T}}_1  = \left[ {x_0 ,\infty } \right) \cap \mathbb{T}_1$ and $\overline{\mathbb{T}}_2 = \left[ {y_0 ,\infty } \right) \cap \mathbb{T}_2$ for $x_0 \in \mathbb{T}_1 $ and $y_0 \in \mathbb{T}_2 $ and $\Omega=\mathbb{T}_1 \times \mathbb{T}_2$.
 Let $\sigma_1,\rho_1,\Delta_1$ and $\sigma_2,\rho_2,\Delta_2$ denote the forward jump operators, backward jump operators and the delta differentiation operator respectively on $\mathbb{T}_1$
 and $\mathbb{T}_2$. Let $a<b$ be points in $\mathbb{T}_1$, $c<d$ are point in  $\mathbb{T}_2$, $[a,b)$ is the half closed bounded interval in $\mathbb{T}_1$, and $[c,d)$ is the half closed bounded interval in $\mathbb{T}_2$.

 We say that a real valued function $f$ on $\mathbb{T}_1\times\mathbb{T}_2$ at $(t_1,t_2) \in \overline{\mathbb{T}}_1 \times \overline{\mathbb{T}}_2 $ has a $\Delta_1$ partial derivative  $f^{\Delta_1}(t_1,t_2)$ with respect to $t_1$ if for each $\epsilon >0$ there exists a neighborhood $U_{t_1 }$ of $t_1$ such that
\[
\left| {f\left( {\sigma _1 \left( {t_1 } \right),t_2 } \right) - f\left( {s,t_2 } \right) - f^{\Delta _1 } \left( {t_1 ,t_2 } \right)\left( {\sigma _1 \left( {t_1 } \right) - s} \right)} \right| \le \varepsilon \left| {\sigma _1 \left( {t_1 } \right) - s} \right|,
\]
for each $s \in U_{t_1}$, $t_2 \in \mathbb{T}_2$.
We say that  $f$ on $\mathbb{T}_1\times\mathbb{T}_2$ at $(t_1,t_2) \in \overline{\mathbb{T}}_1 \times \overline{\mathbb{T}}_2 $ has a $\Delta_2$ partial derivative  $f^{\Delta_2}(t_1,t_2)$ with respect to $t_2$ if for each $\eta >0$ there exists a neighborhood $U_{t_2 }$ of $t_2$ such that
\[
\left| {f\left( {t_1 ,\sigma _2 \left( {t_2 } \right)} \right) - f\left( {t_1 ,l} \right) - f^{\Delta _2 } \left( {t_1 ,t_2 } \right)\left( {\sigma _2 \left( {t_2 } \right) - l} \right)} \right| \le \eta \left| {\sigma _2 \left( {t_2 } \right) - l} \right|,
\]
for all $l \in U_{t_2}$, $t_1 \in \mathbb{T}_1$.
 The function $f$ is called rd-continuous in $t_2$ if for every $\alpha_1 \in \mathbb{T}_1$, the function $f(\alpha_1,.)$ is rd-continuous on $\mathbb{T}_2$. The function $f$ is called rd-continuous in $t_1$ if for every $\alpha_2 \in \mathbb{T}_2$ the function $f(.,\alpha_2)$ is rd-continuous on $\mathbb{T}_1$.

Let $CC_{rd}$ denote the set of functions $f(t_1,t_2)$ on $\mathbb{T}_1 \times \mathbb{T}_2$ where $f$ is rd continuous in $t_1$ and $t_2$.
    Let $CC'_{rd}$ deontes the set of all functions $CC_{rd}$  for which both the $\Delta _1$ partial derivative and $\Delta _2$ partial derivative exists and are in $CC_{rd}$.

	The basic information on time scales and inequalities  can be found in \cite{Agr1,Ana1,Boh1,Boh2}.

\section{\textbf{Ostrowski Inequalities for double integrals on time scales }  }
Now we give Ostrowski Inequalities for double integrals on time scales

\paragraph{\textbf{Theorem 2.1}}
Let  $f,g \in CC'_{rd} \left( {\left[ {a,b} \right] \times \left[ {c,d} \right],\mathbb{R}} \right)$ and  $f^{\Delta _2 \Delta _1 } \left( {x,y} \right)$, $g^{\Delta _2 \Delta _1 } \left( {x,y} \right)$ exist rd-continuous on $[a,b] \times [c,d]$. Then

\begin{align*}
&\left| {\int\limits_a^b {\int\limits_c^d {\left\{ {\left[ {f\left( {x,y} \right)g\left( {x,y} \right) - } \right.} \right.} } } \right.\frac{1}{2}\left[ {P\left( {f\left( {x,y} \right)} \right)g\left( {x,y} \right)} \right. \\
&\left. {\left. { + P\left( {g\left( {x,y} \right)} \right)f\left( {x,y} \right)} \right]} \right\}\left. {\Delta _2 y\Delta _1 x} \right| \\ 
&\le \frac{1}{8}\left( {b - a} \right)\left( {d - c} \right)\int\limits_a^b {\int\limits_c^d {\left[| {g\left( {x,y} \right)}| \right.} } \left\| {f^{\Delta _2 \Delta _1 } } \right\|_\infty   \\ 
& {\left. { + f\left( {x,y} \right)\left\| {g^{\Delta _2 \Delta _1 } } \right\|} \right]_\infty  \Delta _2 y\Delta _1 x} .
\tag{2.1}
\end{align*}
where 
\begin{align*}
P\left( {f\left( {x,y} \right)} \right)
& = \frac{1}{2}\left[ {f\left( {\sigma _1 \left( s \right),y} \right) + f\left( {x,\sigma _2 \left( t \right)} \right) + f\left( {x,d} \right) + f\left( {b,y} \right)} \right] \\
&- \frac{1}{4}\left[ {f\left( {\sigma _1 \left( s \right),\sigma _2 \left( t \right)} \right) + f\left( {\sigma _1 \left( s \right),d} \right) + f\left( {b,\sigma _2 \left( t \right)} \right) + f\left( {b,d} \right)} \right],
\tag{2.2}
\end{align*}
and 
\begin{align*}
&Q\left( {f^{\Delta _2 \Delta _1 } \left( {x,y} \right)} \right) \\
&= \int\limits_{\sigma _1 \left( s \right)}^b {\int\limits_{\sigma _2 \left( t \right)}^d {\frac{{\partial ^2 f\left( {\eta ,\tau } \right)}}{{\Delta _2 \tau \Delta _1 \eta }}} } \Delta _2 \tau \Delta _1 \eta  - \int\limits_{\sigma _1 \left( s \right)}^b {\int\limits_y^d {\frac{{\partial ^2 f\left( {\eta ,\tau } \right)}}{{\Delta _2 \tau \Delta _1 \eta }}} } \Delta _2 \tau \Delta _1 \eta  \\ 
& - \int\limits_x^b {\int\limits_{\sigma _2 \left( t \right)}^d {\frac{{\partial ^2 f\left( {\eta ,\tau } \right)}}{{\Delta _2 \tau \Delta _1 \eta }}} } \Delta _2 \tau \Delta _1 \eta  + \int\limits_x^b {\int\limits_y^d {\frac{{\partial ^2 f\left( {\eta ,\tau } \right)}}{{\Delta _2 \tau \Delta _1 \eta }}} } \Delta _2 \tau \Delta _1 \eta.
\tag{2.3}
\end{align*}
Similarly $P(g(x,y))$ and $Q\left( {f^{\Delta _2 \Delta _1 } \left( {x,y} \right)} \right) $ are defined similar to $(2.2)$ and $(2.3)$.
	\paragraph{\textbf{Proof.}} From the hypotheses we have for $(x,y) \in [a,b] \times [c,d]$
	
	\begin{align*}
	&\int\limits_{\sigma _1 \left( s \right)}^b {\int\limits_{\sigma _2 \left( t \right)}^d {\frac{{\partial ^2 f\left( {\eta ,\tau } \right)}}{{\Delta _2 \tau \Delta _1 \eta }}} } \Delta _2 \tau \Delta _1 \eta  \\
	& = \int\limits_{\sigma _1 \left( s \right)}^x {\left[ {\left. {\frac{{\partial f\left( {\eta ,\tau } \right)}}{{\Delta _1 \eta }}} \right|_{\sigma _2 \left( t \right)}^y } \right]} \Delta _1 \eta  \\ 
	&= \int\limits_{\sigma _1 \left( s \right)}^x {\left[ {\frac{{\partial f\left( {\eta ,y} \right)}}{{\Delta _1 \eta }} - \frac{{\partial f\left( {\eta ,\sigma _2 \left( t \right)} \right)}}{{\Delta _1 \eta }}} \right]} \Delta _1 \eta  \\ 
	&= \left. {f\left( {\eta ,y} \right)} \right|_{\sigma _1 \left( s \right)}^x  - \left. {f\left( {\eta ,\sigma _2 \left( t \right)} \right)} \right|_{\sigma _1 \left( s \right)}^x  \\
	& = f\left( {x,y} \right) - f\left( {\sigma _1 \left( s \right),y} \right) - f\left( {x,\sigma _2 \left( t \right)} \right) + f\left( {\sigma _1 \left( s \right),\sigma _2 \left( t \right)} \right).
	\tag{2.4}
\end{align*}
	Similarly we have 
	\begin{align*}
	&\int\limits_{\sigma _1 \left( s \right)}^b {\int\limits_y^d {\frac{{\partial ^2 f\left( {\eta ,\tau } \right)}}{{\Delta _2 \tau \Delta _1 \eta }}} } \Delta _2 \tau \Delta _1 \eta \\
	&=  - f\left( {x,y} \right) - f\left( {\sigma _1 \left( s \right),d} \right) + f\left( {x,d} \right) + f\left( {\sigma _1 \left( s \right),y} \right),
\tag{2.5}
\end{align*}
\begin{align*}
	 &\int\limits_x^b {\int\limits_{\sigma _2 \left( t \right)}^d {\frac{{\partial ^2 f\left( {\eta ,\tau } \right)}}{{\Delta _2 \tau \Delta _1 \eta }}} } \Delta _2 \tau \Delta _1 \eta   \\
	&=  - f\left( {x,y} \right) - f\left( {b,\sigma _2 \left( t \right)} \right) + f\left( {x,\sigma _2 \left( t \right)} \right) + f\left( {b,\sigma _2 \left( t \right)} \right), 
\tag{2.6}
\end{align*}
and
\begin{align*}
	&\int\limits_x^b {\int\limits_y^d {\frac{{\partial ^2 f\left( {\eta ,\tau } \right)}}{{\Delta _2 \tau \Delta _1 \eta }}} } \Delta _2 \tau \Delta _1 \eta 
	& = f\left( {x,y} \right) + f\left( {b,d} \right) - f\left( {x,d} \right) - f\left( {b,y} \right).
\tag{2.7}
\end{align*}

Adding above identities we have
\begin{align*}
& 4f\left( {x,y} \right) - 2f\left( {\sigma _1 \left( s \right),y} \right) - 2f\left( {x,\sigma _2 \left( t \right)} \right) - 2f\left( {x,d} \right) - 2f\left( {b,y} \right) \\
&+ f\left( {\sigma _1 \left( s \right),\sigma _2 \left( t \right)} \right) + f\left( {\sigma _1 \left( s \right),d} \right) + f\left( {b,\sigma _2 \left( t \right)} \right) + f\left( {b,d} \right) \\
&= \int\limits_{\sigma _1 \left( s \right)}^b {\int\limits_{\sigma _2 \left( t \right)}^d {\frac{{\partial ^2 f\left( {\eta ,\tau } \right)}}{{\Delta _2 \tau \Delta _1 \eta }}} } \Delta _2 \tau \Delta _1 \eta  - \int\limits_{\sigma _1 \left( s \right)}^b {\int\limits_y^d {\frac{{\partial ^2 f\left( {\eta ,\tau } \right)}}{{\Delta _2 \tau \Delta _1 \eta }}} } \Delta _2 \tau \Delta _1 \eta  \\ 
&- \int\limits_x^b {\int\limits_{\sigma _2 \left( t \right)}^d {\frac{{\partial ^2 f\left( {\eta ,\tau } \right)}}{{\Delta _2 \tau \Delta _1 \eta }}} } \Delta _2 \tau \Delta _1 \eta  + \int\limits_x^b {\int\limits_y^d {\frac{{\partial ^2 f\left( {\eta ,\tau } \right)}}{{\Delta _2 \tau \Delta _1 \eta }}} } \Delta _2 \tau \Delta _1 \eta.
\tag{2.8}
\end{align*}
From $(2.3)$, $(2.4)$ and $(2.8)$ we have
\[
f\left( {x,y} \right) - P\left( {f\left( {x,y} \right)} \right) = \frac{1}{4}Q\left( {f^{\Delta _2 \Delta _1 } \left( {x,y} \right)} \right),
\tag{2.9}\]
for $x,y \in [a,b] \times [c,d]$. \\
Similarly for function $g$ we have
\[
g\left( {x,y} \right) - P\left( {g\left( {x,y} \right)} \right) = \frac{1}{4}Q\left( {g^{\Delta _2 \Delta _1 } \left( {x,y} \right)} \right),
\tag{2.10}\]
for $x,y \in [a,b] \times [c,d]$. 
Multiplying $(2.9)$ and $(2.10)$ by $g(x,y)$ and $f(x,y)$ and adding the resulting identities we get
\begin{align*}
&2f\left( {x,y} \right)g\left( {x,y} \right) - g\left( {x,y} \right)P\left( {f\left( {x,y} \right)} \right) - f(x,y)P\left( {g\left( {x,y} \right)} \right) \\ 
&= \frac{1}{4}g\left( {x,y} \right)Q\left( {f^{\Delta _2 \Delta _1 } \left( {x,y} \right)} \right) + \frac{1}{4}f\left( {x,y} \right)Q\left( {g^{\Delta _2 \Delta _1 } \left( {x,y} \right)} \right).
\tag{2.11}
\end{align*}
Integrating $(2.11)$ over $[a,b] \times [c,d]$ we have
\begin{align*}
&\int\limits_a^b {\int\limits_c^d {\left[ {f\left( {x,y} \right)g\left( {x,y} \right)} \right.} }  \\ 
& \left. { - \frac{1}{2}\left[ {P\left( {f\left( {x,y} \right)} \right)g\left( {x,y} \right) + P\left( {g\left( {x,y} \right)} \right)f\left( {x,y} \right)} \right]} \right]\Delta _2 y\Delta _1 x \\ 
&= \frac{1}{8}\int\limits_a^b {\int\limits_c^d {\left[ {Q\left( {f^{\Delta _2 \Delta _1 } \left( {x,y} \right)} \right)g\left( {x,y} \right) + Q\left( {g^{\Delta _2 \Delta _1 } \left( {x,y} \right)} \right)f\left( {x,y} \right)} \right]} } \Delta _2 y\Delta _1 x.
\tag{2.12}
\end{align*}

From the properties of modulus we have

\begin{align*}
\left| {Q\left( {f^{\Delta _2 \Delta _1 } \left( {x,y} \right)} \right)} \right| 
& \le \int\limits_a^b {\int\limits_c^d {\left| {f^{\Delta _2 \Delta _1 } \left( {t,s} \right)} \right|} } \Delta _2 s\Delta _1 t \\
& \le \left\| {f^{\Delta _2 \Delta _1 } } \right\|_\infty  \left( {b - a} \right)\left( {d - c} \right),
\tag{2.13}
\end{align*}
and
\begin{align*}
\left| {Q\left( {g^{\Delta _2 \Delta _1 } \left( {x,y} \right)} \right)} \right| 
&\le \int\limits_a^b {\int\limits_c^d {\left| {g^{\Delta _2 \Delta _1 } \left( {t,s} \right)} \right|} } \Delta _2 s\Delta _1 t \\ 
&\le \left\| {g^{\Delta _2 \Delta _1 } } \right\|_\infty  \left( {b - a} \right)\left( {d - c} \right).
\tag{2.14}
\end{align*}
From $(2.12),(2.13)$ and $(2.14)$ we have
\begin{align*}
&\left| {\int\limits_a^b {\int\limits_c^d {\left[ {f\left( {x,y} \right)} \right.} } } \right.g\left( {x,y} \right) \\ 
&\left. {\left. { - \frac{1}{2}\left[ {P\left( {f\left( {x,y} \right)} \right)g\left( {x,y} \right) + P\left( {g\left( {x,y} \right)} \right)f\left( {x,y} \right)} \right]} \right]\Delta _2 y\Delta _1 x} \right| \\
&\le \frac{1}{8}\int\limits_a^b {\int\limits_c^d {\left[ {\left| {g\left( {x,y} \right)} \right|} \right.} } \left| {Q\left( {f^{\Delta _2 \Delta _1 } \left( {x,y} \right)} \right)} \right| \\ 
& \left. { + \left| {f\left( {x,y} \right)} \right|\left| {Q\left( {g^{\Delta _2 \Delta _1 } \left( {x,y} \right)} \right)} \right|} \right]\Delta _2 y\Delta _1 x \\ 
&\le \frac{1}{8}\int\limits_a^b {\int\limits_c^d {\left[ {\left| {g\left( {x,y} \right)} \right|} \right.} } \int\limits_a^b {\int\limits_c^d {\left| {f^{\Delta _2 \Delta _1 } \left( {t,s} \right)} \right|} } \Delta _2 s\Delta _1 t \\ 
&\left. { + \left| {f\left( {x,y} \right)} \right|\int\limits_a^b {\int\limits_c^d {\left| {g^{\Delta _2 \Delta _1 } \left( {t,s} \right)} \right|} } \Delta _2 s\Delta _1 t} \right]\Delta _2 y\Delta _1 x \\
& \le \frac{1}{8}\left( {b - a} \right)\left( {d - c} \right)\int\limits_a^b {\int\limits_c^d {\left[ {\left| {g\left( {x,y} \right)} \right|} \right.} } \left\| {f^{\Delta _2 \Delta _1 } } \right\|_\infty   \\ 
&\left. { + \left| {f\left( {x,y} \right)} \right|\left\| {g^{\Delta _2 \Delta _1 } } \right\|_\infty  } \right|\Delta _2 y\Delta _1 x.
\tag{2.15}
\end{align*}
Which is required inequality. \\
\paragraph{}
  Now we give the continuous and Discrete equivalent version  of above inequality where $\mathbb{T}=\mathbb{R}$ and $\mathbb{T}=\mathbb{Z}$
\paragraph{\textbf{Corollary 2.1}}(Continuous Case)
If we put $\mathbb{T}_1=\mathbb{T}_2=\mathbb{R}$ we have 
\begin{align*}
&\left| {\int\limits_a^b {\int\limits_c^d {\left[ {f\left( {x,y} \right)} \right.} } } \right.g\left( {x,y} \right) \\ 
&\left. {\left. { - \frac{1}{2}\left[ {P\left( {f\left( {x,y} \right)} \right)g\left( {x,y} \right) + P\left( {g\left( {x,y} \right)} \right)f\left( {x,y} \right)} \right]} \right]dydx} \right| \\
&\le \frac{1}{8}\left( {b - a} \right)\left( {d - c} \right)\int\limits_a^b {\int\limits_c^d {\left[ {\left| {g\left( {x,y} \right)} \right|} \right.} } \left\| {D_2 D_1 f} \right\|_\infty   \\
&\left. { + \left| {f\left( {x,y} \right)} \right|\left\| {D_2 D_1 g} \right\|_\infty  } \right|dydx. 
\end{align*}

\begin{align*}
P(f(x,y)) 
&= \frac{1}{2}\left[ {f\left( {x,c} \right) + f\left( {x,d} \right) + f\left( {a,y} \right) + f\left( {b,y} \right)} \right] \\ 
&- \frac{1}{4}\left[ {f\left( {a,c} \right) + f\left( {a,d} \right) + f\left( {b,c} \right) + f\left( {b,d} \right)} \right],
\end{align*}

\begin{align*}
Q\left( {D_2 D_1 f\left( {x,y} \right)} \right)  
&= \int\limits_a^x {\int\limits_c^y {D_2 D_1 f\left( {t,s} \right)} } dsdt - \int\limits_a^x {\int\limits_y^d {D_2 D_1 f\left( {t,s} \right)} } dsdt \\
&- \int\limits_x^b {\int\limits_c^y {D_2 D_1 f\left( {t,s} \right)} } dsdt + \int\limits_x^b {\int\limits_y^d {D_2 D_1 f\left( {t,s} \right)} } dsdt .
\end{align*}

Similarly $P(g(x,y))$ and $Q(D_2 D_1 g(x,y))$,
which is Ostrowski inequality for Double integral.

\paragraph{\textbf{Corollary 2.2}}(Discrete Case)
If $\mathbb{T}_1=\mathbb{T}_2=\mathbb{Z}$ and $a=c=0$, $b=k \in \mathbb{N}$ and $d=r \in \mathbb{N}$. Then
\begin{align*}
&\left| {\sum\limits_{x = 1}^k {\sum\limits_{y = 1}^r {\left\{ {f\left( {x,y} \right)g\left( {x,y} \right)} \right.} } } \right. 
 \left. { - \frac{1}{2}\left( {g\left( {x,y} \right)P\left( {f\left( {x,y} \right)} \right) + f\left( {x,y} \right)P\left( {g\left( {x,y} \right)} \right)} \right)} \right| \\
&\le \frac{1}{8}kr\sum\limits_{x = 1}^k {\sum\limits_{y = 1}^r {\left[ {\left| {g\left( {x,y} \right)} \right|\left\| {\Delta _2 \Delta _1 f} \right\|_\infty   + \left| {f\left( {x,y} \right)} \right|\left\| {\Delta _2 \Delta _1 g} \right\|_\infty  } \right]} }., 
\end{align*}

\begin{align*}
P\left( {f\left( {x,y} \right)} \right)
&= \frac{1}{2}\left[ {f\left( {x,1} \right) + f\left( {x,r + 1} \right) + f(1,y) + f(k + 1,y)} \right] \\ 
&- \frac{1}{4}\left[ {f\left( {1,1} \right) + f\left( {1,r + 1} \right) + f(k + 1,1) + f(k + 1,r + 1)} \right],
\end{align*}

\begin{align*}
Q\left( {\Delta _2 \Delta _1 f\left( {x,y} \right)} \right)
&\sum\limits_{s = 1}^{x - 1} {\sum\limits_{t = 1}^{y - 1} {\Delta _2 \Delta _1 f\left( {s,t} \right)} }  - \sum\limits_{s = 1}^{x - 1} {\sum\limits_{t = 1}^m {\Delta _2 \Delta _1 f\left( {s,t} \right)} }  \\
&- \sum\limits_{s = x}^k {\sum\limits_{t = 1}^{y - 1} {\Delta _2 \Delta _1 f\left( {s,t} \right)} }  + \sum\limits_{s = x}^k {\sum\limits_{t = y}^r {\Delta _2 \Delta _1 f\left( {s,t} \right)} }.
\end{align*}
Which is Discrete Ostrowski Inequality. \\
 Now we give some Ostrowski inequality for double integrals.
  \paragraph{\textbf{Theorem 2.2}}
Let $f,g,P(f(x,y)),P(g(x,y),f^{\Delta _2 \Delta _1 },g^{\Delta _2 \Delta _1 }$ be as in Theorem $2.1$ then
\begin{align*}
\left| {\int\limits_a^b {\int\limits_c^d {\left\{ {f\left( {x,y} \right)g\left( {x,y} \right)} \right.} } } \right.
&- \left[ P \right.\left( {f\left( {x,y} \right)} \right)g\left( {x,y} \right) + P\left( {g\left( {x,y} \right)} \right)f\left( {x,y} \right) \\ 
&\left. {\left. { - P\left( {f\left( {x,y} \right)} \right)P\left( {g\left( {x,y} \right)} \right)} \right\}\Delta _2 y\Delta _1 x} \right| \\ 
&\le \frac{1}{{16}}\left\{ {\left( {b - a} \right)\left( {d - c} \right)} \right\}^2 \left\| {f^{\Delta _2 \Delta _1 } } \right\|_\infty  \left\| {g^{\Delta _2 \Delta _1 } } \right\|_\infty .
\tag{2.16}
\end{align*}
for $[x,y] \in [a,b] \times [c,d]$.
\paragraph{\textbf{Proof.}}
From $(2.9)$ and $(2.10)$ we have
\[
f(x,y) - P\left( {f\left( {x,y} \right)} \right) = \frac{1}{4}Q\left( {f^{\Delta _2 \Delta _1 } \left( {x,y} \right)} \right),
\tag{2.17}
\]
and
\[
g(x,y) - P\left( {g\left( {x,y} \right)} \right) = \frac{1}{4}Q\left( {g^{\Delta _2 \Delta _1 } \left( {x,y} \right)} \right),
\tag{2.18}
\]
for $[x,y] \in [a,b] \times [c,d]$.
\paragraph{}
Multiplying left hand side and right hand side of $(2.17)$ and $(2.18)$ we get
\begin{align*}
f(x,y)g(x,y)
&- \left[ {f(x,y)P\left( {g\left( {x,y} \right)} \right) + g(x,y)P\left( {f\left( {x,y} \right)} \right)} \right] \\ 
&= \frac{1}{{16}}Q\left( {f^{\Delta _2 \Delta _1 } \left( {x,y} \right)} \right)Q\left( {g^{\Delta _2 \Delta _1 } \left( {x,y} \right)} \right).
\tag{2.19}
\end{align*}
Integrating $(2.19)$ over $[a,b] \times [c,d]$ and from the properties of modulus we have
\begin{align*}
&\left| {\int\limits_a^b {\int\limits_c^d {\left\{ {f\left( {x,y} \right)g\left( {x,y} \right)} \right.} } } \right. \\
& - \left[ P \right.\left( {f\left( {x,y} \right)} \right)g\left( {x,y} \right) + P\left( {g\left( {x,y} \right)} \right)f\left( {x,y} \right) \\ 
&\left. {\left. { - P\left( {f\left( {x,y} \right)} \right)P\left( {g\left( {x,y} \right)} \right)} \right\}\Delta _2 y\Delta _1 x} \right| \\
&\le \frac{1}{{16}}\int\limits_a^b {\int\limits_c^d {\left| {Q\left( {f^{\Delta _2 \Delta _1 } \left( {x,y} \right)} \right)} \right|} } \left| {Q\left( {g^{\Delta _2 \Delta _1 } \left( {x,y} \right)} \right)} \right|\Delta _2 y\Delta _1 x.
\tag{2.20}
\end{align*}

Now using $(2.13)$ and $(2.14)$ in $(2.15)$ we get the required inequality $(2.16)$.
\paragraph{}Now we give continuous and discrete version of the inequality $(2.16)$ where $\mathbb{T}=\mathbb{R}$ and $\mathbb{T}=\mathbb{Z}$ which is as follows
\paragraph{\textbf{Corollary 2.3}} (Continuous Case)
If we put $\mathbb{T}_1=\mathbb{T}_2=\mathbb{R}$ in above we get
\begin{align*}
&\left| {\int\limits_a^b {\int\limits_c^d {\left\{ {f\left( {x,y} \right)g\left( {x,y} \right)} \right.} } } \right. \\
&- \left[ P \right.\left( {f\left( {x,y} \right)} \right)g\left( {x,y} \right) + P\left( {g\left( {x,y} \right)} \right)f\left( {x,y} \right) \\
&\left. {\left. { - P\left( {f\left( {x,y} \right)} \right)P\left( {g\left( {x,y} \right)} \right)} \right\}dydx} \right| \\ 
 & \le \frac{1}{{16}}\int\limits_a^b {\int\limits_c^d {\left\| {D_2 D_1 f} \right\|_\infty  \left\| {D_2 D_1 g} \right\|_\infty  } } dydx. 
\end{align*}
where $f,g,P,Q,D_2 D_1 f,D_2 D_1 g$ is as in Corollary $2.1$.
Which is Ostrowksi type inequality for double integral. 
\paragraph{\textbf{Corollary 2.4}} (Discrete Case)
If $\mathbb{T}_1=\mathbb{T}_2=\mathbb{Z}$ and $a=c=0$,$b=k \in N$ and $d=r \in N$. Then
\begin{align*}
 \left| {\sum\limits_{x = 1}^k {\sum\limits_{y = 1}^r {\left\{ {f\left( {x,y} \right)g\left( {x,y} \right)} \right.} } } \right.
& - \left[ {g\left( {x,y} \right)P\left( {f\left( {x,y} \right)} \right) + f\left( {x,y} \right)P\left( {g\left( {x,y} \right)} \right)} \right. \\ 
&\left. {\left. {\left. { - P\left( {f\left( {x,y} \right)} \right)P\left( {g\left( {x,y} \right)} \right)} \right\}} \right]} \right| \\
&\le \frac{1}{{16}}\left( {kr} \right)^2 \int\limits_a^b {\int\limits_c^d {\left\| {\Delta _2 \Delta _1 f} \right\|_\infty  \left\| {\Delta _2 \Delta _1 g} \right\|_\infty  } }.
\end{align*}
where $P,Q$ are as in Corollary $2.2$. Which is discrete Ostrowski type inequality.

\section{Trapezoid type Inequality on time scales}
Now we give the dynamic Trapezoid type inequality on time scales.
\paragraph{\textbf{Theorem 3.1}}
Let $f,f^{\Delta _2 \Delta _1 }$ be as in Theorem $2.1$. Then
\begin{align*}
&\left| {\int\limits_a^b {\int\limits_c^d {f\left( {t,s} \right)\Delta _2 s\Delta _1 } } } \right.t - \frac{1}{2}\left[ {\left( {d - c} \right)\int\limits_a^b {f\left( {t,\sigma _2 \left( t \right)} \right) + f\left( {t,d} \right)} } \right] \\ 
&\left. { + \left( {b - a} \right)\int\limits_c^d {\left[ {f\left( {\sigma _1 \left( s \right),s} \right) + f\left( {b,s} \right)} \right]} \Delta _2 s} \right] \\ 
&+ \frac{1}{4}\left( {b - a} \right)\left( {d - c} \right)\left[ {f\left( {\sigma _1 \left( s \right),\sigma _2 \left( t \right)} \right) + f\left( {\sigma _1 \left( s \right),d} \right)} \right. \\ 
&\left. { + f\left( {b,\sigma _2 \left( t \right)} \right) + f\left( {b,d} \right)} \right] \\ 
&= \frac{1}{4}\int\limits_a^b {\int\limits_c^d {\left| {f^{\Delta _2 \Delta _1 } \left( {t,s} \right)} \right|} } \Delta _2 s\Delta _2 t. 
\tag{3.1}
\end{align*}

\paragraph{\textbf{Proof.}}
 From the proof of Theorem $2.1$ we have
\begin{align*}
f\left( {x,y} \right)
&= \frac{1}{2}\left[ {f\left( {\sigma _1 \left( s \right),y} \right) + f\left( {x,\sigma _2 \left( t \right)} \right) + f\left( {x,d} \right) + f\left( {b,y} \right)} \right] \\ 
&+ \frac{1}{4}\left[ {f\left( {\sigma _1 \left( s \right),\sigma _2 \left( t \right)} \right) + f\left( {\sigma _1 \left( s \right),d} \right) + f\left( {b,\sigma _2 \left( t \right)} \right) + f\left( {b,d} \right)} \right] \\
&= \frac{1}{4}Q\left( {f^{\Delta _2 \Delta _1 } \left( {x,y} \right)} \right),
 \tag{3.2}
\end{align*}
for $[x,y] \in [a,b] \times [c,d]$.\\
Integrating $(3.2)$ over $[a,b] \times [c,d]$ we get
 \begin{align*}
&\int\limits_a^b {\int\limits_c^d {f\left( {t,s} \right)} } \Delta _2 s\Delta _1 t \\ 
& - \frac{1}{2}\left[ {\left( {d - c} \right)\int\limits_a^b {\left[ {f\left( {t,\sigma _2 \left( t \right)} \right) + f\left( {t,d} \right)} \right]} } \right. \Delta _1 t\\ 
& + \left. {\left( {b - a} \right)\int\limits_c^d {\left[ {f\left( {\sigma _1 \left( s \right),s} \right) + f\left( {b,s} \right)} \right]\Delta _2 s} } \right] \\ 
&+ \frac{1}{4}\left( {b - a} \right)\left( {d - c} \right)\left[ {f\left( {\sigma _1 \left( s \right),\sigma _2 \left( t \right)} \right) + f\left( {\sigma _1 \left( s \right),d} \right)} \right. \\
&\left. { + f\left( {b,\sigma _2 \left( t \right)} \right) + f\left( {b,d} \right)} \right] \\ 
& = \frac{1}{4}\int\limits_a^b {\int\limits_c^d {Q\left( {f^{\Delta _2 \Delta _1 } \left( {t,s} \right)} \right)} } \Delta _2 s\Delta _1 t.
\tag{3.3}
\end{align*}
From the property of modulus and integrals we have
\[
\left| {Q\left( {f^{\Delta _2 \Delta _1 } \left( {x,y} \right)} \right)} \right| \le \int\limits_a^b {\int\limits_c^d {\left| {f^{\Delta _2 \Delta _1 } \left( {t,s} \right)} \right|} } \Delta _2 s\Delta _1 t.
\tag{3.4}\]
From $(3.3)$ and $(3.4)$ we have

 \begin{align*}
&\left| {\int\limits_a^b {\int\limits_c^d {f\left( {t,s} \right)} } \Delta _2 s\Delta _1 t} \right. \\
& - \frac{1}{2}\left[ {\left( {d - c} \right)\int\limits_a^b {\left[ {f\left( {t,\sigma _2 \left( t \right)} \right) + f\left( {t,d} \right)} \right]} } \right.\Delta _1 t \\ 
& + \left. {\left( {b - a} \right)\int\limits_c^d {\left[ {f\left( {\sigma _1 \left( s \right),s} \right) + f\left( {b,s} \right)} \right]\Delta _2 s} } \right] \\ 
&+ \frac{1}{4}\left( {b - a} \right)\left( {d - c} \right)\left[ {f\left( {\sigma _1 \left( s \right),\sigma _2 \left( t \right)} \right) + f\left( {\sigma _1 \left( s \right),d} \right)} \right. \\ 
&\left. {\left. { + f\left( {b,\sigma _2 \left( t \right)} \right) + f\left( {b,d} \right)} \right]} \right| \\ 
&  \le \frac{1}{4}\int\limits_a^b {\int\limits_c^d {\left| {Q\left( {f^{\Delta _2 \Delta _1 } \left( {t,s} \right)} \right)} \right|} } \Delta _2 s\Delta _1 t \\ 
& \le \frac{1}{4}\left( {b - a} \right)\left( {d - c} \right)\int\limits_a^b {\int\limits_c^d {\left| {f^{\Delta _2 \Delta _1 } \left( {t,s} \right)} \right|} } \Delta _2 s\Delta _1 t.
\tag{3.5}
\end{align*}
which is required inequality.
Now we give the Continuous and discrete version of Trapezoid inequality when  $\mathbb{T}=\mathbb{R}$ and $\mathbb{T}=\mathbb{Z}$.
\paragraph{\textbf{Corollary 3.1}} (Continuous Case)
If we put $\mathbb{T}_1=\mathbb{T}_2=\mathbb{R}$ in above we get
\begin{align*}
& \int\limits_a^b {\int\limits_c^d {f\left( {t,s} \right)dsdt} } 
- \frac{1}{2}\left[ {\left( {d - c} \right)\int\limits_a^b {\left[ {f\left( {t,c} \right) + f\left( {t,d} \right)} \right]} dt} \right. \\ 
&\left. { + \left( {b - a} \right)\int\limits_c^d {\left[ {f\left( {a,s} \right) + f\left( {b,s} \right)} \right]ds} } \right] \\ 
&+ \frac{1}{4}\left( {b - a} \right)\left( {d - c} \right)\left[ {f\left( {a,c} \right) + f\left( {a,d} \right) + f\left( {b,c} \right) + f\left( {b,d} \right)} \right] \\ 
& \le \frac{1}{4}\left( {b - a} \right)\left( {d - c} \right)\int\limits_a^b {\int\limits_c^d {\left| {D_2 D_1 f\left( {t,s} \right)} \right|} } dtds.
\end{align*}
which is Continuous Trapezoid type inequality

\paragraph{\textbf{Corollary 3.2}} (Discrete Case)
If $\mathbb{T}_1=\mathbb{T}_2=\mathbb{Z}$ then we get
\begin{align*}
&\left| {\sum\limits_{t = 1}^k {\sum\limits_{s = 1}^r {f\left( {t,s} \right)} } } \right. - \frac{1}{2}\left[ {m\sum\limits_{t = 1}^k {\left[ {f\left( {t,1} \right) + f\left( {t,m + 1} \right)} \right]} } \right. \\ 
&\left. { + k\sum\limits_{s = 1}^r {\left[ {f\left( {1,s} \right) + f\left( {k + 1,s} \right)} \right]} } \right] \\
&+ \frac{1}{4}km\left[ {f\left( {1,1} \right) + f\left( {1,r + 1} \right) + f\left( {k + 1,1} \right) + f\left( {k + 1,r + 1} \right)} \right].
\end{align*}
which is Discrete Trapezoid type inequality

\paragraph{\textbf{Acknowledgment}}
This research is supported by Science and Engineering Research Board (SERB, New
Delhi, India) Project File No. SB/S4/MS:861/13.

\end{document}